\documentclass[10pt]{article}

\usepackage{amsmath}
\usepackage{amssymb}
\usepackage{latexsym}

\usepackage[all]{xy}

\newtheorem{thm}{Theorem}[section]
\newtheorem{prop}[thm]{Proposition} \newtheorem{lemma}[thm]{Lemma}
\newtheorem{cor}[thm]{Corollary} \newtheorem{dfn}[thm]{Definition}
 \newtheorem{rmk}[thm]{Remark}
\newtheorem{ex}[thm]{Example} 

\newcommand {\pf}{\noindent{\bf Proof.}\ }
\newcommand{\complex}{{\mathbb C}}

\newcommand{\reals}{{\mathbb R}}

\newcommand{\Hom}{{\rm Hom}}

\newcommand{\backl}{\mathbin{\vrule width1.5ex height.4pt\vrule height1.5ex}}

\newcommand{\cale}{{\cal E}}

\newcommand{\cali}{{\cal I}}
\newcommand{\calj}{{\cal J}}

\newcommand{\calt}{{\cal T}}

\newcommand{\calu}{{\cal U}}
\newcommand{\calw}{{\cal W}}
\newcommand{\del}{\partial}
\newcommand{\delbar}{\overline{\partial}}
\newcommand{\qed}{\begin{flushright} $\Box$\ \ \ \ \ \end{flushright}}

\newcommand{\arrows}{\,\lower1pt\hbox{$\longrightarrow$}\hskip-.24in\raise2pt
             \hbox{$\longrightarrow$}\,}

\newcommand{\mlp}{-\log\psi}

\newcommand{\vbar}{\overline{v}}

\newcommand{\gambar}{\overline{\gamma}}

\begin{document}

\title{{\bf Poisson geometry\\
    and deformation quantization \\near a strictly pseudoconvex boundary}}
\author
{Eric Leichtnam
\\ CNRS and Institut Math\'ematique de Jussieu\\
Etage 7 E,
175 Rue du Chevaleret\\
Paris 75013  France
\\ 
\\
Xiang Tang
\\ Department of Mathematics\\
University of California, Davis, CA, 95616 USA
\\{\small (xtang@math.ucdavis.edu)}
\\
\\
 Alan Weinstein
\thanks{Research partially supported by NSF Grant
DMS-0204100
\newline \mbox{~~~~}MSC2000 Subject Classification Numbers: 32T15 (Primary); 53D10,  53D17
\newline \mbox{~~~~}Keywords: Poisson structure, pseudoconvexity,
plurisubharmonic function, contact structure, Lie algebroid}\\
Department of Mathematics\\ University of California, Berkeley, CA
94720 USA\\ {\small(alanw@math.berkeley.edu)}}
\date{}
\maketitle



\vspace{-.25in}

\begin{abstract}
Let $X$ be a complex manifold with strongly pseudoconvex boundary $M$.
If $\psi$ is a defining function for $M$, then $-\log\psi$ is
plurisubharmonic on a neighborhood of $M$ in $X$, and the
(real) 2-form $\sigma = i \del \delbar(-\log \psi)$ is a symplectic
structure on the complement of $M$ in a neighborhood in $X$ of $M$;
it blows up along $M$.

The Poisson structure obtained by
inverting $\sigma$ extends smoothly across $M$ and determines
a contact structure on $M$ which is the same as the one induced by the
complex structure.  When $M$ is compact,
the Poisson structure near $M$
is completely determined up to isomorphism
by the contact structure on $M$.
In addition, when $-\log\psi$ is plurisubharmonic
throughout $X$, and $X$ is compact, bidifferential
operators constructed by Engli\v{s} for the
Berezin-Toeplitz deformation quantization of $X$ are smooth up to the
boundary.  The proofs use a complex Lie algebroid determined by the CR
structure on $M$, along with some ideas of Epstein, Melrose, and
Mendoza concerning manifolds with contact boundary.

\end{abstract}
\section{Introduction}
\label{sec-intro}
Let $u$ be a
real-valued function on a complex manifold $X$.  The  2-form
$\sigma=\sigma_u=i\del\delbar u$ is of type (1,1), real, and exact
(since $\del\delbar = d \delbar$).  In addition,
$\sigma(x,y)=\sigma(Jx,Jy)$, where $J$ is the complex structure
viewed  as an endomorphism of $TX$, so the ``hermitian hessian''
bilinear form  $g_u(x,y)=\sigma(x, J y)$ is symmetric.
When $g_u$ is  positive definite, $u$ is said to be {\bf
strongly plurisubharmonic}.  In this  case, $g_u$ is a K\"ahler metric on $X$,
and the form $\sigma$ is  nondegenerate, i.e. symplectic.  The
function $u$ is called a {\bf K\"ahler  potential} for $g_u$.
Conversely, the Dolbeault  lemma implies that any K\"ahler metric on
$X$ arises from a potential  on a neighborhood of each point of $X$.
While the local geometry of the   K\"ahler metric depends very much on
the choice of potential, the symplectic form has no local invariants.
 (See \cite{el-gr:convex} for global results about the symplectic
geometry of K\"ahler manifolds with global potential functions.)

In this paper, we will investigate what happens when the complex manifold $X$ has
a boundary $M$.  The hyperplane field
$F_{M,X} = TM \cap JTM$ is called the {\bf maximal complex subbundle} of $TM$.
A {\bf defining function} for $M$
 is a smooth nonnegative
function $\psi$ whose zero set is $M$, and which has no critical points on
$M$.  The restriction to $F_{M,X}$ of the hermitian hessian $g_{\psi}$
 is, up to a
positive conformal factor, independent of the choice of defining
function.  The associated invariant object, a symmetric bilinear form
with values in the conormal bundle of $M$, is called the {\bf Levi
form} of $M$.  $M$ is called {\bf Levi nondegenerate}
when this form is nondegenerate and {\bf strongly pseudoconvex} when
it is negative definite.  $M$ is Levi nondegenerate if and only if
$F_{M,X}$ is a contact structure.

When $M$ is strongly pseudoconvex, $-\log\psi$ is strongly
 plurisubharmonic on $U\setminus M$ for some neighborhood $U$ of $M$
 in $X$, and
 $\sigma_{\mlp}$ on $U\setminus M$ is a symplectic structure
which blows up along $M$.
We will show that the corresponding Poisson
 structure $\pi_{\mlp}$ extends smoothly to $M$, along which it is zero.

In fact, using
Epstein, Melrose and Mendoza's \cite{ep-me-me:resolvent}
notion of $\Theta$ structure, we will show that, whenever $M$ is Levi nondegenerate,
the
local isomorphism type of the Poisson structure $\pi_{\mlp}$ is
independent of everything but the dimension of $X$.  The local model
is Lebrun's \cite{le:thickenings} Poisson structure
on the normal
bundle to a contact structure, and equivalence with this model gives
the smoothness of $\pi_{\mlp}$ up to the boundary.

Contact structures also play a role in
the global (on $M$) version of this result:
the germ along $M$ of the
Poisson structure is determined, up to diffeomorphisms fixing $M$,
by the contact structure
$F_{M,X}$.
Its isomorphism class is thus
independent of the choice of the defining function and of the choice
of (compatible) complex structure.  Again, the model for $\pi_{\mlp}$ is given by
LeBrun's construction.

Similar results, for flows and infinitesimal deformations on
pseudoconvex manifolds,
have been obtained by Kor\'anyi and Reimann \cite
{ko-re:contact}\cite{re:quasiconformal}.

We turn next to quantization.  When $X$ is compact and $\psi$ is
strongly plurisubharmonic throughout $X$, $\pi_{\mlp}$ is the
semiclassical commutator of the Berezin-Toeplitz deformation
quantization product on $X$.  The construction of this product 
involves the action of smooth 
functions on $X$ by multiplication and projection on a
parameterized family of weighted Bergman spaces of holomorphic
functions on the interior of $X$.  The Berezin-Toeplitz
 product was analyzed in
the pseudoconvex setting by Engli\v{s} \cite{en:berezin},
following many earlier studies on closed manifolds.  He showed
that the induced product on smooth functions has an asymptotic
expansion in the weight parameter; the terms in the expansion
are bidifferential
operators whose coefficients are algebraic combinations of
the K\"ahler metric, its curvature, and covariant derivatives
thereof.

We will use the notion of {\bf complex Lie
algebroid} to show that all the bidifferential operators in the
Berezin-Toeplitz-Engli\v{s} quantization are smooth up to the
boundary.   In fact, we will show something stronger.   Karabegov
\cite{ka:separation} has defined a notion of quantization
with separation of variables on a K\"ahler manifold,
and it is known that the Berezin-Engli\v{s}-Toeplitz quantization has this property
on the interior of $X$.  We will extend Karabegov's definition by
by introducing a notion of 
{\bf para-K\"ahler structure} on a complex Lie algebroid, and
we will show that the separation of variables property 
then holds up to the boundary.

In the future, we plan to extend our results to more general
manifolds, replacing the global strongly plurisubharmonic function
$\psi$ by a family of local functions obtained from a connection on a
hermitian line bundle.  These results will be used in the proof of a
topological formula for the relative index of CR structures defined by
Epstein \cite{ep:relative} which was conjectured by Atiyah and
Weinstein \cite{we:questions}.  (A proof of the conjecture by Epstein
\cite{ep:subelliptic} has recently appeared, but his methods are quite
different from ours.)

\section{Normal forms}

For most of   this section, we will forget about complex geometry and
look at what Epstein, Melrose, and Mendoza \cite{ep-me-me:resolvent}
call $\Theta$-structures.
$M$ will now be the boundary of any manifold $X$ of real dimension $2n+2.$

\begin{dfn}
\label{dfn-EMM}
An {\bf EMM form} is a 1-form $\Theta$ on $X$ whose pullback to $M$ is a
contact form.
\end{dfn}

Our model example of an EMM form will be the pullback to $X=M\times
\reals^+$ of a contact form $\phi$ on $M$; we identify $M$ with the
zero set of the coordinate function $r$
 on $\reals^+ =  [0,\infty)$.  The symplectic form
$d(\Theta/r)$  blows up
along $M$, but
LeBrun \cite{le:thickenings} observed that
the Poisson structure inverse to $d(\Theta/r)$
extends to a smooth Poisson structure on $M\times \reals^+$.
If $\phi = du+\sum p_jdq^j$ in local coordinates $(u,q,p)$ on $M$,
then, on $M\times\reals^+$,
\begin{equation}
\label{eq-colebrun}
d(\Theta/r)=(1/r^2)\left[-dr\wedge\left(du+\sum p_j dq^j\right)+r\-dp_j\wedge dq^j\right],
\end{equation}
and the Poisson structure corresponding\footnote{There is a choice of
  sign when one says that a Poisson structure corresponds to a
  symplectic structure.  Unlike LeBrun \cite{le:thickenings}, we use
  the convention in which $dq\wedge dp$ corresponds to the relation
  $\{q,p\}=1$.} to $-d(\Theta/r)$ is
\begin{equation}
\label{eq-lebrun}
\Pi=r\left[\left(r\frac{\partial}{\partial
r} + \sum p_{j}\frac{\partial}{\partial p_{j}}\right) \wedge
\frac{\partial}{\partial u} + \sum \frac{\partial}{\partial q^{j}}
\wedge \frac{\partial}{\partial p_{j}}\right].
\end{equation}

\subsection{Local normal form}
We will use the following local theorem in order
to obtain a global normal form.  (It would be nice to get
the global form all at once, but we do not know how to do it.)

\begin{thm}
\label{thm-local}
Let $\Theta$ be an EMM form on the manifold $X$ with boundary $M$, and
let $\psi$ be a defining function for the boundary.
Then, near each $m\in M$,
 there exist local coordinates
$(q,p,u,r)$ on $X$
in which $\psi=r$ and
$d(\Theta/\psi)$ has the form \eqref{eq-colebrun}
 on the
complement of $M$.  In
particular, $-d(\Theta/\psi)$ is symplectic on the complement of $M$ in a
neighborhood of $m$ in $X$; the corresponding
Poisson structure on this neighborhood has the local normal form
\eqref{eq-lebrun}.
\end{thm}

 \pf
 We begin by setting $r=\psi$.
 By the Darboux theorem for
contact 1-forms \cite{ca:lectures}, we may find local coordinates on
$M$ for which the pullback of $\Theta$ has the expression
 $du+p_j
dq^j.$ In fact, the pullback of
 $\Theta$ to each level of $r$ near
$M$ is still a contact structure, so
 we may choose coordinates on
all these levels, depending smoothly on
 $r$, so that the pullbacks
all have the same form.  It follows that
 $\Theta$ itself may be
written as $du+p_j dq^j + a dr,$ where $a$ is a
 smooth function of
all the variables.

To eliminate the term $a dr$, we use Moser's method, i.e. constructing
a diffeomorphism (preserving $r$ and fixed on $M$) by integrating
a time-dependent vector field $X_t$.  As usual, we define $\Theta_t$
by interpolation as $du+p_j dq^j+ t a \- dr$ and choose $X_t$ to
satisfy the condition $X_t \backl d(\Theta_t/r) = -(a/r)dr. $
Now
\begin{align}
d(\Theta_t/r)=&(1/r)d\Theta_t  -(1/r^2)dr\wedge \Theta_t\\
             =&(1/r)(dp_j\wedge dq^j + t da\wedge dr)-(1/r^2)dr\wedge
(du+p_j dq^j+ t a \-  dr). \notag
\end{align}
If we take $X_t $ to be a function $ f_t(q,p,u,r)$
times the (Reeb)  vector field $\del/\del u,$ it
will be tangent to the levels of $r$.  $X_t$ must satisfy the equation
$$X_t \backl  d(\Theta_t/r) =
(t/r)(f_t\, \del a/\del u)
+(1/r)f_t)dr = -(a/r)dr,$$
which has the solution
$f_t=
-ra/(1+rt \,\del a/\del u).$  The denominator is invertible near $M$, and
the factor of $r$ in the numerator of $f_t$ insures that $X_t$
vanishes along $M$, in addition to being smooth and tangent to
the levels of $r$.
\qed

\subsection{Global normal form}

To put $\Theta$ in normal form on a neighborhood of the
 entire boundary, we can no longer fix the $\psi$ levels, because the
 characteristic line element field of the pullback of $d(\Theta/\psi)$ has a
 global dynamics which may vary from one $\psi$ level to another.
 This also makes it impossible to use the Darboux theorem as we did
 for the local normal form.  Instead, we use Gray's theorem, which
 asserts that deformations of a contact {\em structure} on a compact
 manifold are trivial.  We refer the reader to Cannas da Silva
\cite{ca:lectures} for a proof, noting
 for use below that the transformations in Gray's theorem may easily be
 chosen to depend smoothly on a parameter.

Like the local model, the global normal form comes from LeBrun
\cite{le:thickenings}.  For any contact structure  $F\subset TM$, the conormal
bundle $\nu^* = (TM/F)^*$ may be identified with the 1-dimensional subbundle
of $T^*M$ consisting of all real multiples of any contact form
defining the contact structure.  The pullback to $\nu^*$ of the
canonical symplectic structure on $T^*M$ is nondegenerate on the complement of
the zero section of $\nu^*$.  We may identify this complement by
``inversion'' with the complement of the zero section in the normal
bundle $\nu=TM/F$.
LeBrun shows that the Poisson structure corresponding to this form
on $\nu$ now extends smoothly over the zero section.  When $\nu$ is
oriented, a
choice of contact form identifies $\nu$ with $M\times \reals$, and
the nonnegative normal bundle $\nu^+$ is identified with $M\times
\reals^+$.  The Poisson structure is given by \eqref{eq-lebrun}.

\begin{thm}
\label{thm-global} Let $\Theta$ be an EMM form for the manifold
$X$ with boundary $M$, and $\psi$ a defining function for the
boundary. Then  a neighborhood of $M$ in $X$, with the structure
corresponding to $-d(\Theta/\psi)$, is Poisson isomorphic to  a
neighborhood of the zero section in the nonnegative normal bundle
$\nu^+$, with
 the LeBrun-Poisson structure associated to the contact structure
 induced by $\theta$ on $M$.
\end{thm}

\pf
As in the proof of Theorem \ref{thm-local}, we will work with the
2-forms.  Since the diffeomorphism we construct will be smooth along
the singular locus of these forms, it will automatically be a Poisson
isomorphism.

To begin, we identify both $X$ near $M$ and $\nu^+$ near the zero
section with a neighborhood of the zero section in the trivial bundle
$M\times \reals^+$.  For
$\nu^+$,  we use the trivialization of $\nu$ given by the contact
form which is the pullback of $\Theta$ to $M$.
For $X$ near $M$,  we first arrange that the projection onto
$\reals$ is the given function $\psi$, i.e. we set $r=\psi$;
we then use Gray's theorem to
arrange that the projection of each $\psi$ level onto $M$ is a contact
diffeomorphism.

Let us write $\Theta_0$ for the standard form $du+\sum p_jdq^j$ (independent of $r$)
 and $\Theta_1$ for the
given form.  Since $\Theta_1$ defines the same contact structure as
$\Theta_0$ on each level of $r$, and it agrees with $\Theta_0$ on the zero level,
$\Theta_1=\Theta_0+a\,dr +br\Theta_0$, where $a$ and $b$ are smooth
functions.
As before, we linearly interpolate to
get $\Theta_t=\Theta_0+ta\,dr +tbr\Theta_0$.  Note for later use in this
proof
that these are all EMM
forms, so we can apply Theorem \ref{thm-local} to
put them in local normal form.

Once again, we seek a time-dependent vector field $X_t$ to generate
our normalizing transformation.  The required condition on this
vector field is
$$X_t\backl d(\Theta_t/r) = - (a/r)dr + b\Theta_0.$$
The unique solution of this equation is (on the complement of
$M$) the contraction of the
right hand side with the Poisson structure corresponding to
$-d(\Theta_t /r)$.
We already know from the local normal form that this Poisson structure
vanishes along $M$, so the contraction extends smoothly over $M$.  It
remains to show that the contraction vanishes along $M$.  For this, it
suffices to show that the contraction with $dr$ vanishes to second
order.  But, from the local normal form
\eqref{eq-lebrun}, we find immediately that, in normal form
coordinates, this contraction is equal to $r^2 du$, and our proof is complete.
\qed

\begin{rmk}
\label{rmk-rigid} {\em Although the normal form theorem above may
suggest that the LeBrun-Poisson structure is rigid with respect to
arbitrary higher-order perturbations, this is in fact not the
case.  For instance, when $M$ is 1-dimensional, the Poisson
structure is simply $\Pi=r^2 \frac{\del}{\del r}\wedge
\frac{\del}{\del u}$.  This structure is exact in the sense that
there is a vector field $\xi$ (namely $\frac{\del}{\del r}$)
satisfying $[\xi,\Pi]=\Pi,$ but the Poisson structure
$(r^2+r^3)\frac{\del}{\del r}\wedge \frac{\del}{\del u}$ does not
admit such a $\xi$ if $M$ is a circle.  A related fact is that
$(1/(r^2+r^3))dr\wedge du$ is not $d(\Theta /r)$ for any EMM form
$\Theta$. }
\end{rmk}

We also note the following relative form of Theorem \ref{thm-global}.

\begin{cor}
\label{cor-jet}
If $\Theta_0$ and $\Theta_1$ are EMM forms which agree to infinite order along $M$, then there is a diffeomorphism germ on $X$ along $M$ which agrees with the identity to infinite order along $M$ and pulls back $\Theta_1$ to $\Theta_0$.
\end{cor}

\pf
By Theorem \ref{thm-global}, we may assume that $X = M \times [0,1)$ and that $\Theta_0$ is the standard form.  We then repeat the proof of Theorem \ref{thm-global}; the functions $a$ and $b$ vanish to infinite order along $M$, hence so does the vector field $X_t$ which generates the normalizing transformation.
\qed

\subsection{Application to pseudoconvex boundaries}

Let $X$ be a complex manifold with boundary $M$.
By elementary calculus on complex manifolds,
$$\delbar f = \frac{1}{2} (df + i J^* df)$$
for any smooth function $f$, where $J:TX\to TX$ is the almost complex
structure.  It follows that
$$ i \del\delbar (\mlp) = i d\delbar(\mlp)=\frac{1}{2}d(J^*\frac{
  d\psi}{\psi}) =
d (\frac{1}{2}J^*d\psi/\psi).  $$
Let $\psi$ be a defining function for the boundary and set $r=\psi$
  and $\Theta=\frac{1}{2}J^*d\psi$.  We will show that $\Theta$ is an
  EMM form.  In fact, in $T_M X$, $TM$ is $\ker d \psi$, so
$$\ker \Theta \cap TM  = \ker J^*d\psi\cap TM=J(\ker d\psi)\cap
  TM=JTM\cap TM=F_{M,X},$$ the maximal complex subbundle of $TM$.
  Since $M$ is Levi-nondegenerate, $F_{M,X}$ is a contact structure,
  and hence $\Theta$ is an EMM form.
It follows that all the results of this section apply to the form
  $\sigma_{\mlp}$ and the corresponding Poisson structure
  $\pi_{\mlp}$.

We remark that our results correspond very closely to results on flows
and deformations due to Kor\'anyi and Reimann \cite{ko-re:contact}
\cite{re:quasiconformal}.  Since the Poisson structure determines the
contact structure on the boundary, our methods also give a simple
proof of their (easier) converse result that a smooth map which is
symplectic on the interior must be contact on the boundary.

Here is a direct nondegeneracy proof which is
independent of the normal form theorem.  It involves a volume element
computation which we will use in the proof of Proposition
\ref{prop-canonical} below.

\begin{prop}
\label{prop-nondegen}
Let $\psi$ be a defining function for the boundary $M$ of $X$.  The closed 2-form
$\sigma_{\mlp}$ is nondegenerate on a neighborhood of $M$ in the
interior of $X$ if $M$ is Levi-nondegenerate.
\end{prop}

\pf
We compute:
$$(1/i)\sigma_{\mlp}=\del\delbar( \mlp) = -\del\left( \frac{\delbar
  \psi}{\psi}\right) = \frac{-\psi\del\delbar \psi +
  \del\psi\wedge\delbar\psi}{\psi^2}.$$

Raising this 2-form to the $n+1$'st power gives
$$\left(-\frac{\del\delbar\psi}{\psi} +
\frac{\del\psi\wedge\delbar\psi}{\psi^2}\right)^{n+1}
=\left(-\frac{\del\delbar\psi}{\psi}\right)^{n+1}+
\left(-\frac{\del\delbar\psi}{\psi}\right)^{n}\wedge
\frac{\del\psi\wedge\delbar\psi}{\psi^2}
$$
which is $\psi^{-(n+2)}$ times
$$\psi\left(-\del\delbar\psi\right)^{n+1}+\left(-\del\delbar\psi\right)^{n}\wedge\del\psi\wedge\delbar\psi.$$

Our lemma will be proven if we can show that this form is
nonzero near $M$ when $M$ is Levi-nondegenerate.
Since the first term vanishes along $M$ it suffices to show that the
second is nonzero.
Since  $\del\psi\wedge\delbar\psi$ is nonzero and annihilates
the contact structure $F_{M,X}$, the nonvanishing of the term
is equivalent to nondegeneracy
of the restriction to $F_{M,X}$ of
$\left(-\del\delbar\psi\right)^{n}$.  But the latter is just the
2-form associated via $J$ to the Levi form.
\qed

\begin{rmk}
\label{rmk-bundles}

{\em There are natural 1-1 correspondences among several bundles along
  the boundary $M$ whose sections admit natural simple and transitive
  actions of the smooth positive functions on $M$.

\begin{enumerate}
\item
1-jets along $M$ of defining functions.
\item
Sections of the conormal bundle $TM^\perp\subset T^*_M X$ which are
``positive'' in the sense that they take positive values on
inward-pointing vectors in $T_M X.$
\item
Contact forms realizing the cooriented contact structure on $M$.
\item
Volume elements on $M$ compatible with the natural boundary
orientation of $M$.
\end{enumerate}

The correspondence $1\leftrightarrow 2$ is almost tautological, since
any section of the conormal bundle may be realized as the derivative
along $M$ of a defining function. (For instance, if we multiply the
defining function $\psi$ by a positive function $\lambda$, its
differential along $M$ is also multiplied by $\lambda$.)  For
$2\leftrightarrow 3$, we associate to each positive section $\alpha$
of the conormal bundle the pullback to $M$ of $J^* \alpha$.  (To go in
the other direction, we extend any contact form along $M$ to a section
of $T^*_M X$ by requiring it to annihilate $JTM$.)
Finally, for $3\leftrightarrow 4$, we associate to each contact form
$\theta$ the volume element $\theta\wedge(d\theta)^n.$  Replacing
$\theta$ by $\lambda$ multiplies the volume element by
$\lambda^{n+1}$.
}
\end{rmk}

\section{Geometry on complex Lie algebroids}

Complex Lie algebroids were defined in \cite{ca-we:geometric}
and have been studied in more detail in \cite{bl:duality} and \cite{we:complex}.  In this section,
we will review the definitions and use a complex Lie algebroid
to ``regularize'' the geometry of a complex manifold near a
pseudoconvex boundary.

\subsection{Definition and first examples}

We recall that a {\bf Lie algebroid} over a smooth manifold $X$ is a
real vector
bundle $E$ over $X$ with a Lie algebra
structure (over $\reals$) on its sections and with a bundle map
$\rho$ (called the {\bf anchor}) from $E$ to the tangent
bundle $TM$, satisfying
the Leibniz rule
$$[a,fb]=f[a,b]+(\rho(a)f)b$$ for sections $a$ and $b$ and smooth functions
$f$.

There is an analogous definition for complex manifolds, in which $E$
is a holomorphic vector bundle over $X$, and the Lie algebra structure
is defined on the sheaf of local sections.   Such objects are called
complex Lie algebroids by Chemla \cite{ch:duality}, but, as in
\cite{ca-we:geometric}, we will reserve this term
for the ``hybrid'' concept defined below.

\bigskip

   From now on, $C^\infty(X)$ will denote the algebra of smooth
{\em complex-valued} functions on a manifold $X$.

\begin{dfn}
A {\bf complex Lie algebroid} over a smooth (real) manifold $X$ is a complex
vector bundle $E$ over $X$ with a Lie algebra structure (over $\complex$)
on its space $\cale$ of sections and a bundle map $\rho$ (called the {\bf
anchor}) from $E$ to the complexified tangent bundle $T_\complex X$, satisfying
 the Leibniz rule
$$[a,fb]=f[a,b]+(\rho(a)f)b$$ for $a$ and $b$ in $\cale$ and $f$
in $C^\infty(X)$.
\end{dfn}

A ``trivial'' class of complex Lie algebroids consists of the
complexifications of real Lie algebroids, such as  $T_\complex X$
itself.  More interesting are general ``involutive systems,''
which are subbundles of  $T_\complex X$ whose spaces of
sections are closed under the
(complexified) bracket of vector fields.  (Up to isomorphism, these are just the complex Lie algebroids with injective anchor.)  Among these are the
complex structures and CR structures.  By a complex structure, we mean here
a subbundles of the form $E=T^{0,1}_J X = \{v+iJv|v\in TX\}$, where
$J:TX\to TX$ is an integrable almost complex structure.
These are characterized among all complex subbundles
by closure under bracket and the algebraic property that  $T_\complex M= E
\oplus \overline E.$  By a CR structure, we mean an involutive system
$E$ for which $E\cap\overline E=\{0\}$ and $E + \overline E$ has
codimension $1$ in $T_\complex X$.  Any real hypersurface $M$ in a
complex manifold $X$ (such as a boundary)
inherits a CR structure, namely the intersection
$G_{M,X}=T_\complex M \cap T^{0,1}_J X$.  (The problem of realizing a
given CR structure in this way has been crucial in the development of
linear PDE theory.)  The sum  $G_{M,X}\oplus \overline{G_{M,X}}$
is the complexification of the maximal complex subbundle
 $F_{M,X}$.

The main example of our paper, introduced in Section
\ref{subsec-natural}, will not have an injective anchor.  However, its
anchor will be bijective on an open dense subset of the base manifold
$X$, and the use of Lie algebroids with this property could be viewed
as an application of the method of moving frames, extended to allow
certain ``singular" frame fields.

No discussion of complex Lie algebroids should fail to mention the
important example of generalized complex structures
\cite{gu:generalized}, \cite{hi:generalized}, but having thus
fulfilled this obligation, we will not discuss them further.

\subsection{Some constructions on complex Lie algebroids}
Many notions can be extended from real to complex Lie algebroids
without any extra effort.  Here are some
which we will use later.   Parts of this section are
almost transcribed verbatim from \cite{ne-ts:quantization}.
Note that  all the constructions below are local and may thus be
carried out on the sheaf level.

\begin{dfn}
\label{dfn-ederham}Let $(E, \rho, [,])$ be a complex Lie algebroid
over $X$. The $E$-{\bf de Rham complex} $(^E\Omega ^\bullet(X), ^E d)$ is given
by $^E
\Omega^\bullet(X)=\Gamma(\wedge^\bullet(E^*))$, with
$$
^E d\mu(a_1, \ldots, a_{k+1})=\sum_i
(-1)^i\rho(a_i) \mu(a_1, \ldots, \hat{a}_i, \ldots,
a_{k+1})+$$
$$\sum_{i<j}(-1)^{i+j-1}\mu([a_i, a_j],
a_1, \ldots, \hat{a}_i, \ldots, \hat{a}_j,\ldots,
a_{k+1}).$$  Elements of the complex are called $E$-differential forms on $X$;
the cohomology of $^E d$ is denoted by $^E H^\bullet(X)$ and is called the
$E$-de Rham cohomology of $X$.
\end{dfn}

When $E=T_\complex X$, $^E H^\bullet(X)$ is the usual de Rham cohomology of
$X$ with complex coefficients.

\begin{dfn}
An $E$-connection on a vector bundle $F$ over $X$ is a map
$$(a,\gamma)\mapsto \nabla_a \gamma$$ from $\cale \times\Gamma(F)$ to
$\Gamma(F)$ which is $C^\infty (X)$-linear in $a$ and  satisfies
the Leibniz rule $$\nabla_a f\gamma = f\nabla_a \gamma +
(\rho(a)f)\gamma$$ for $f\in C^\infty(X)$. Equivalently, an
$E$-connection on  $F$  is a map $\nabla: \Gamma (F) \rightarrow
\Gamma( E^* \otimes F)$ satisfying $\nabla ( f \gamma) = f \nabla
\gamma + ^{E} d f  \otimes \gamma.$

\end{dfn}

Like an ordinary linear connection, an $E$-connection extends to a map
\[
\nabla:\Gamma(\wedge^\bullet(E^*)\otimes F)\to
\Gamma(\wedge^{*+1}(E^*)\otimes F).
\]
The square of this extended operator is given by
$
\nabla^2 \gamma = R\wedge \gamma,
$
where the {\bf curvature} $R$ is the element of $\wedge^2(E^*)\otimes
\mathrm{End}(F)$
 defined by
\[
R(a,b)=\nabla_a\nabla_b-\nabla_b\nabla_a-\nabla_{[a, b]}.
\]
When the
curvature is zero, the connection is also called a {\bf
  representation} of $E$ on $F$.

When $F=E$, we may also define the {\bf torsion} of $\nabla$ by the
formula:
$T(a,b)=\nabla_a b - \nabla_b a - [a,b].$  As in the case of the
tangent bundle, the torsion is a skew-symmetric tensor, i.e. a section of
$\wedge^2(E^*)\otimes E.$  The usual construction of the Levi-Civita
connection applies, so that, given a field of nondegenerate symmetric
inner products on $E$, there is a unique connection without torsion
which is compatible with the inner product.

\begin{ex}
\label{ex-bott}
{\em The flat ``Bott connection'' on the normal bundle to a foliation is the
linearization of the holonomy.  But the construction is purely formal
and can be extended to the situation where $E'$ is any subalgebroid of
a Lie algebroid $E$.  Namely, we define an $E'$-connection on the
quotient vector bundle $E/E'$ by the rule $\nabla_a<b>=<[a,b]>$, where
$a$ and $b$ are sections of $E$ and $<\cdot>$ denotes the equivalence
class modulo $E'$.  (We use angled instead of the usual square brackets for the
equivalence class to avoid
confusion with the Lie algebroid operation.)  It is straightforward to
check that this $\nabla$ is a Lie algebroid representation.

We note that $\nabla$ can be seen as the representation on
homology, as in Appendix A of \cite{ev-lu-we:transverse},
associated to a natural representation up to homotopy of $E$ on
the short complex of Lie algebroids $0\to E'\to E\to 0$
}
\end{ex}

We also introduce  a universal enveloping
object first defined in a slightly different way by Rinehart \cite{ri:differential}.
\begin{dfn}
\label{dfn-eop}Let $(E, \rho, [,])$ be a complex Lie algebroid on
$X$, and $\calt$  the free associative (i.e. tensor) algebra with
generators in
$C^{\infty}(X)$ (of degree 0) and
$\Gamma (E)$ (of degree 1). The algebra
$^E Op$  of $E$-differential operators on $X$ is defined as
$\calt/\cali$,
where $\cali$ is the two-sided
ideal of $\calt$ generated by elements of the
form
$$f\otimes g-fg,$$
$$f\otimes
a-fa,$$
$$a\otimes b-b\otimes a-[a, b],$$
  and
$$a\otimes (fb)-(fa)\otimes b-(\rho(a))f)b,$$
for $a, b \in \cale$ and $f, g \in C^\infty(X)$.
\end{dfn}

The grading of $\calt$ defines a filtration $^E Op_n$ of $^E
Op$, and the following result is a straightforward application of Theorem 3.1
in \cite{ri:differential}.
\begin{lemma}
\label{lem:pbw}For any complex Lie algebroid $(E, \rho, [,])$ over $X$,
there is a natural isomorphism $Gr^EOp(X)\simeq \Gamma(X, S(E))$,
where $S(E)$ is the bundle of symmetric algebras on the fibres of
$E$.  In particular, the algebra $C^\infty(X)$ may be identified with a
subalgebra of $^EOp$.
\end{lemma}

Following Calaque \cite{ca:formality}, we may also introduce the space
of $E$-polydifferential (or multi-differential) operators with its
Gerstenhaber bracket operation.  (This structure was already suggested by
Xu \cite{xu:quantum} and used implicitly by Nest and Tsygan \cite{ne-ts:quantization}.)

The usual jet spaces of functions on $X$ are not sensitive enough to the
action of $E$, since sections of the isotropy act trivially, 
so we must use the following generalization.

\begin{dfn}
\label{dfn-ejets} Let $(E, \rho, [, ])$ be a complex Lie algebroid
over $X$. The space of $E$-jets on $X$ is the linear space
$^E\!{\calj}(X)=\Hom_{C^{\infty}(X)}(^E Op(X), C^{\infty}(X))$.
\end{dfn}

In the real case, the $E$-jets may be identified with the jets of
functions along the units of a (local) groupoid integrating $E$.
A similar identification also works in the complex case, though the
integration in the sense of \cite{we:complex} may be only formal.  

The complex analog of Proposition 2.7 in \cite{ne-ts:quantization} is:

\begin{prop}
\label{prop-flatconn} $^E\calj(X)$ introduced in Definition
\ref{dfn-ejets} is the space of global sections of a
profinite-dimensional vector bundle $^E Jets$.
\end{prop}

We define
the ``Grothendieck connection''
$\nabla^G:\cale \times ^E\!\!\calj(X)\to ^E\!\!\calj(X)$ by
\[
(\nabla_G(a)(l))(D)=a(l(D))-l(a(D)),
\]
for $l\in \Gamma(^E Jets)(=^E\!\!\calj(X) ),$  $a \in \cale$ and $D \in ^E Op(X).$
 As in the real case, this is
a flat connection.

\subsection{A natural complex Lie algebroid on a complex manifold with boundary}
\label{subsec-natural}
We introduce here the complex Lie algebroid which will be central in what follows.   Let $X$ be a complex manifold of dimension $n+1$
 with boundary $M$, and
let $\cale _{M,X}$ be the space of complex vector
fields on $X$ (i.e. sections of $T_\complex X$) whose values along $M$
lie in the induced CR structure $G_{M,X}$. $\cale _{M,X}$ is a module
over $C^{\infty}(X)$  and is closed under bracket.
lemma shows that $\cale _{M,X}$ may be identified with the space
with the space of sections of a complex Lie algebroid $E_{M,X}$.

\begin{lemma}
\label{lemma-module}
$\cale _{M,X}$ is a locally free $C^\infty(X)$-module.
\end{lemma}
\pf Away from the boundary, $\cale _{M,X}$ is the same as
$T_\complex M$, hence locally free. Near a boundary point, we may
choose a local basis $\vbar_1,\ldots,\vbar_n$ of $G_{M,X}$, which
we then extend to a linearly independent set of sections of
$T^{0,1}X$, still denoted by $\vbar_j$, defined in an open subset
of $X$.  We leave the name of the open subset unspecified and will
shrink it as necessary. Let $v_j$ be the complex conjugate of
$\vbar_j$.  These vectors all annihilate $\psi$ on $M$; there is
no obstruction to having them annihilate $\psi$ everywhere. Next,
we choose a local section $\vbar_0$ of $T^{0,1}X$ such that
$\vbar_0\cdot \psi = 1$,  and we let $v_0$ be its conjugate.  This
gives a local basis $(v,\vbar)$ for the complex vector fields.
Such a vector field belongs to $\cale_{M,X}$ if and only if, when
it is expanded with respect to this basis, the coefficients of
$\vbar_0$ and all the $v_j$ vanish along $M$.  Since this means
that all these coefficients are divisible by $\psi$ with smooth
quotient, we get a local basis $(u,u')$ for $\cale_{M,X}$ by
setting $u'_0=\psi\vbar_0$, $u'_j=\vbar_j$ for $j=1,\ldots,n$, and
$u_j=\psi v_j$ for $j=0,\ldots,n$. \qed

The local basis $(u,u')$ constructed in the proof above
may be thought of as a moving frame, some of whose entries vanish
along $M$.   The crucial property here is that
 the structure functions which express
Lie brackets  in the given frame are smooth up to $M$.

We note that the complex conjugates of the basis vectors are
$\overline{u_0}=u'_0$ and
$\overline{u_j}=\psi u'_j$ for $j=1,\ldots,n$.  The Lie algebroid
$E_{M,X}$ does {\em not} admit an operation of complex conjugation.

We will also use the coframe $(\theta,\theta')$ dual to $(u,u')$.
Denoting by $(\gamma,\gambar)$ the basis of complex-valued 1-forms dual
  to $(v,\vbar)$, we find that  $\gamma^0=\del\psi$ and
 $\gambar^0=\delbar \psi$.  For the vector bundle $E_{M,X}^*$ dual
to $E_{M,X}$,
we get the local basis of sections
$\theta^j=(1/\psi)\gamma^j$ for $j=0,\ldots,n$ (so that
$\theta^0=\del(\log|\psi|), $
$\theta'^0=(1/\psi)\gambar^0=(1/\psi)\delbar\psi=\delbar(\log|\psi|)$, and
$\theta'^j = \gambar^j$ for $j=1,\ldots,n$.   The complex conjugates are
$\overline{\theta^0}=\theta'^0$ and
$\overline{\theta^j}=(1/\psi)\theta'^j$ for $j=1,\ldots,n$.

The coframe $(\theta,\theta')$ is an ordinary coframe on the interior
of $X$.  Some of these forms blow up along $M$, but the structure
functions which express the exterior differentials of these forms in
terms of the coframe are smooth up to $M$.

\subsection{Para-K\"ahler Lie algebroids}
Recall that a pseudo-K\"ahler structure on a manifold $X$ is a
symplectic structure together with a totally complex polarization.
This means that we have a (real) nondegenerate closed 2-form $\omega$
on $X$ and an integrable subbundle $T^{0,1}X$ of $T_\complex X$ which
is isotropic with respect to the complex extension of $\omega$ and for
which $T_\complex X = \overline{T^{0,1}X}\oplus T^{0,1}X$.  We write
$T^{1,0}X$ for $\overline{T^{0,1}X}.$  The structure is K\"ahler when the
nondegenerate quadratic form $g$ defined on $T^{1,0}X$ by $g(u,v) =
\omega (u,J\overline{v})$   is positive definite.

Thinking of a general complex Lie algebroid $E$ over $X$ as a
substitute for $T_\complex X$, it is natural to try to define an
analogous notion of K\"ahler structure, but we lack the
operation of complex conjugation.  On the other hand, the study of
pseudo-K\"ahler structures often makes little or no use of the
quadratic form $g$, but only of the nondegenerate pairing between
$T^{1,0}X$ and $T^{0,1}X$ defined by the restriction of the symplectic
form.  We are thus dealing with a generalization (by complexification
and passing from tangent bundles to general Lie algebroids) of the
so-called ``para-K\"ahler'' \cite{li:probleme} or ``bilagrangian''
\cite{he:connections} structures, which consist of a
symplectic form together with a transverse pair of lagrangian foliations.

The following definitions are useful in both the real and complex cases.

\begin{dfn}
\label{dfn-parakahler}
A {\bf [complex] symplectic Lie algebroid} is a [complex] Lie algebroid
 $E$ together with an $E$-differential 2-form $\omega$ which is $^E d$
 closed and nondegenerate.  A {\bf polarization} of $(E,\omega)$ is a
 lagrangian subalgebroid of $E$, i.e. a subbundle which is closed
 under brackets and  maximal isotropic with respect
 to $\omega$.  A {\bf [complex] para-K\"ahler Lie algebroid} is a [complex]
 symplectic Lie algebroid with a splitting $E=E^{1,0}\oplus E^{0,1}$
 as the direct sum of two polarizations.

\end{dfn}

\begin{rmk}
\label{rmk-pairing}
{\em
The restriction of $\omega$ to $E^{1,0}\times E^{0,1}$ is a
nondegenerate pairing which we will continue to denote by $\omega$.
This pairing is  also the restriction of a unique {\em symmetric} inner
product on $E$ for which $E^{1,0}$ and $E^{0,1}$ are isotropic.
}
\end{rmk}

\begin{ex}
\label{ex-splitting}
{\em The complex Lie algebroid $E_{M,X}$ of Section \ref{subsec-natural} is naturally split as a direct sum $E_{M,X}^{1,0}\oplus E_{M,X}^{0,1}.$
Sections of $E_{M,X}^{1,0}$ are fields of holomorphic tangent vectors
which vanish on the boundary, while sections of $E_{M,X}^{0,1}$ are
fields  of antiholomorphic tangent vectors which are tangent to the
boundary.\footnote{In the language of \cite{me:geometric}, we are
  dealing with a hybrid of the $0$-calculus and the $b$-calculus.}
Near the boundary, the sections of the two summands
are spanned by the $\theta^j$ and $\theta'^j$ respectively, for
$j=0,\ldots,n$.

In Section \ref{sec-quasik}, we will construct a complex symplectic
structure for which these summands become lagrangian.  }
\end{ex}

\subsection{The para-K\"ahler connection}
\label{sec:para-conn}
The complexification of the Levi-Civita connection on a
pseudo-K\"ahler manifold has many nice properties with respect to
the splitting of the complexified tangent bundle into its holomorphic
and antiholomorphic summands.  In fact, it can be constructed directly
from this splitting and from the pairing given by the complexified
symplectic structure.  By imitating this construction, we may
construct on any para-K\"ahler Lie algebroid $E$ a torsion free
$E$-connection which is compatible with the para-K\"ahler
structure.  (There is in fact just one connection with these properties.)
In the case of a bilagrangian manifold, the
construction yields the bilagrangian connection of Hess \cite{he:connections},
and in fact, beyond a change of terminology, there is nothing we do
here which is not taken from this special case.

\begin{prop}
\label{prop-connection} Let $(E=E^{1,0}\oplus E^{0,1}, \omega)$ be
a para-K\"ahler Lie algebroid.  There is a unique torsion-free
$E$-connection $\nabla$ on $E$ for which covariant differentiation
leaves the para-K\"ahler structure invariant; i.e. for any $a,b,c
\in \Gamma(E)$, $\nabla_a$ leaves the splitting invariant, and
$\rho(a)(\omega(b,c)) =\omega(\nabla_a b ,c)+\omega(b,\nabla_a
c).$  The curvature of this connection $\nabla$ is in
$({E^{1,0}}^*\wedge {E^{0,1}}^*)\otimes {\rm End }(E)$.
\end{prop}

\pf
Our connection $\nabla$ will be built from two flat partial
connections on $E$ defined on the summands.

First, identifying
$E^{0,1}$ with $E/E^{1,0}$, we have via
Example \ref{ex-bott} an $E^{1,0}$-connection on $E^{0,1}$.  Writing
$p^{1,0}$ and $p^{0,1}$ for the projection maps associated to the
splitting of $E$, we therefore have
$$\nabla_a b' = p^{0,1}[a,b']$$ for $a\in \Gamma(E^{1,0})$ and $b' \in
\Gamma(E^{0,1})$.

This partial connection induces a connection on the dual bundle to
$E^{0,1}$, which we identify with $E^{1,0}$ via the pairing
$\omega$.  The resulting $E^{1,0}$-connection on $E^{1,0}$
is determined by the equation
$$\omega(\nabla_a b,c') = \rho(a)(\omega(b,c') - \omega(b,[a,c']).$$
We recall that $\rho$ is the anchor of the Lie algebroid and that we
may omit the projection from the last term because $E^{1,0}$ is
isotropic for $\omega$.

Putting together these two pieces, we get an $E^{1,0}$ connection
$\nabla^{1,0}$ on $E$ which is clearly compatible with the
para-K\"ahler structure.

Now we may interchange the two summands and repeat everything
above to get the required $E^{0,1}$ connection $\nabla^{0,1}$ on
$E$, and then assemble everything to get the required
$E$-connection on $E$.  We leave to the reader the exercise of
verifying (using the fact that $\omega$ is a {\em closed} 2-form),
that this connection has zero torsion.

By the Jacobi identity for $E^{1,0}$ and $E^{0,1}$ vector fields,
we find that the curvatures of the connections $\nabla^{1,0}$ and
$\nabla^{0,1}$ in
$(\wedge^2({E^{1,0}}^*)+\wedge^2({E^{0,1}}^*))\otimes {\rm
End}(E)$ vanish. This implies that the curvature $R$ of the
$E$-connection $\nabla=\nabla^{1,0}+\nabla^{0,1}$ is a $(1,1)$
form. \qed

We call this connection the {\bf para-K\"ahler connection}.

\begin{rmk}
\label{rmk-levicivita}
{\em
A remark made in \cite{et-sa:canonical} is still valid here: the
para-K\"ahler connection is the Levi-Civita connection of the
symmetric inner product of Remark \ref{rmk-pairing}.  This must be so,
since the Levi-Civita connection is unique, and the symmetric inner
product, being built in a canonical way from the para-K\"ahler
structure, must be invariant under the para-K\"ahler connection.
}
\end{rmk}

Given bases $\xi_0,\ldots,\xi_n$ and $\xi'_0,\ldots,\xi'_n$
of $E^{1,0}$ and $E^{0,1}$ respectively, we will write
$\omega_{ij}=\omega(\xi_i , \xi'_j)$ and $\pi^{ij}$ for the
inverse matrix.  We may expand the brackets between summands as
$[\xi_i,\xi'_j] = m_{ij}^k \xi_k + {m'}_{ij}^k \xi'_k$.
It is straightforward to derive the following formulas for the
para-K\"ahler  connection.

$$\nabla_{\xi_i} \xi_j =
\pi^{lk}(\rho(\xi_i)\omega_{jl}-\omega_{jp}{m'}_{il}^p)\xi'_k$$
$$\nabla_{\xi_i} \xi'_j = {m'}_{ij}^k \xi_k.$$

In the usual pseudo-K\"ahler case, we can choose $\xi_i=\partial/\partial
z^i$ and $\xi'_i=\partial/\partial \overline{z}^i$ to
make all the brackets vanish, in which case we get the familiar formulas

$$\nabla_{\xi_i} \xi_j = \pi^{lk}(\partial \omega_{jl}/\partial z^i)
\xi_k.$$
$$\nabla_{\xi_i} \xi'_j = 0.$$

On the other hand, in any para-K\"ahler Lie algebroid, we may choose
the bases 
$\xi_0,\ldots,\xi_n$ and $\xi'_0,\ldots,\xi'_n$ to be dual to one
another with respect to the pairing, so that $\omega_{ij}$ and $\pi^{ij}$ are
identity matrices.  Then we get
\[
\begin{split}
\nabla_{\xi_i }\xi_j &= \sum\limits_{k}-{m'}_{ik}^j \xi_k,\\
\nabla_{\xi_i} \xi'_j&= \sum\limits_k {m'}_{ik}^j \xi'_k.
\end{split}
\]

As a result of the above observation, we see that all the
calculations of K\"ahler geometry can be carried out in an
arbitrary para-K\"ahler Lie algebroid.  This has the following
consequence.

\begin{cor}
\label{cor-para}
Let $X$ be a manifold (possibly with boundary), $E\to X$ a
para-K\"ahler Lie algebroid whose anchor $\rho:E\to T_\complex X$ is
invertible on an open dense subset $\calu \subset X$.  Suppose that
the induced para-K\"ahler structure on $T_\complex \calu$ comes from a
pseudo-K\"ahler structure on $\calu$.  Then all contravariant tensors
and multi-differential operators on $\calu$ which are constructed from
the complex structure, the pseudo-K\"ahler metric, its curvature and
covariant derivatives thereof are the image under $\rho$ of smooth
objects defined on all of $E$.  In particular, they extend smoothly
from $\calu$ to $X$.
\end{cor}

\begin{rmk}
\label{rmk-smoothness}
{\em
One may apply this corollary to the Berezin transform and all the coefficients in the Berezin and Berezin-Toeplitz products,
as analyzed by Engli\v s.  This proves smoothness up
to the boundary of these constructions without any extra work.
Also, we can get smoothness of the canonical form (see Proposition
\ref{prop-canonical}), since it is the Ricci form of the canonical
connection. However, to identify the Berezin-Toeplitz product with
an $E$-product in Theorem \ref{thm:smooth-boundary}, we need the
machinery of formal integrals used by Karabegov and Schlichenmaier
\cite{ka-sch:identification}.
}
\end{rmk}

\begin{rmk}
\label{rmk-analytic}
{\em
It is a much more delicate problem to decide, when the pseudo-K\"ahler
structure on $\calu$ is positive definite, whether elliptic analysis
can be used as in the compact K\"ahler setting to get results 
valid on all of $X$.  For instance, under what conditions
on the singularities of $\rho$ does the deformation quantization with
separation of variables, which extends smoothly to $X$,
arise from a Berezin-Toeplitz symbol calculus as is the case
in Theorem \ref{thm:smooth-boundary} below?
}
\end{rmk}

\subsection{The para-K\"ahler Lie algebroid near a pseudoconvex boundary}
\label{sec-quasik}

We return to our study of the Poisson geometry of a
complex manifold $X$ with pseudoconvex boundary $X$, now using the complex Lie
algebroid $E_{M,X}$.  We begin by showing that the symplectic
structure $\sigma_{\mlp}$, which is singular along $M$, is
perfectly regular as an $E_{M,X}$-symplectic structure.

\begin{thm}
\label{thm-algebroid}
The pullback of $\sigma_{\mlp}$ to
$E_{M,X}$ is a smooth section of $\bigwedge^2E_{M,X}^*.$
This section is
nondegenerate along $M$, hence on a neighborhood of $M$, if and only if $M$ is
Levi-nondegenerate.
\end{thm}
\pf
We will express $\sigma_{\mlp}$ in terms of our bases of
sections of $T_\complex^*X$ and $E_{M,X}^*$.  First of all,
with sums over repeated indices ranging from $1$ to $n$, we have
$$\del\delbar\psi = a\gamma^0\wedge\gambar^0 +
b_k\gamma^0\wedge\gambar^k + \overline{b}_k\gamma^k\wedge\gambar^0 +
c_{jk}\gamma^j\wedge\gambar^k,$$
where $a$ is real and the matrix $c_{jk}$ is hermitian; it is the matrix of
the Levi form.
   From this we get
$$(1/i)\sigma_{\mlp}=-(1/\psi)( a\gamma^0\wedge\gambar^0 + b_k
\gamma^0\wedge\gambar^k + \overline{b}_k\gamma^k\wedge\gambar^0 +
c_{jk}\gamma^j\wedge\gambar^k) + (1/\psi^2)\gamma^0\wedge\gambar^0$$
$$=-(1/\psi)(a\psi\theta^0\wedge\psi\theta'^0
+b_k\psi\theta^0\wedge \theta'^k + \overline{b}_k
\psi\theta^k\wedge\psi\theta'^0 +
c_{jk}\psi\theta^j\wedge\theta'^k) +
(1/\psi^2)\psi\theta^0\wedge\psi\theta'^0$$
$$=(1-\psi a)\theta^0\wedge\theta'^0-b_k\theta^0\wedge \theta'^k
+ \overline{b}_k\psi\theta^k\wedge\theta'^0-c_{jk}\theta^j\wedge\theta'^k.$$
Along $M$, where $\psi=0$, this becomes
$$\theta^0\wedge(\theta'^0-b_k\theta'^k)-  c_{jk}\theta^j\wedge\theta'^k.$$
This is clearly smooth as a section of $\bigwedge^2E_{M,X}^*$, and its
nondegeneracy is  equivalent to that of the
matrix $c_{jk}$, i.e. to that of the Levi form.
\qed

We get another proof of the smooth extension theorem obtained earlier from normal form theory.

\begin{cor}
\label{cor-poisson}
The Poisson structure $\pi_{\mlp}$ obtained by inverting
$\sigma_{\mlp}$ near $M$ extends smoothly over $M$.
\end{cor}
\pf
On the complement of $M$, near $M$, this Poisson structure is the push-forward
of the section of $\bigwedge^2E_{M,X}$ obtained by inverting the
pullback of $\omega_{\mlp}$.  We have just seen that the inverse of
this pullback extends smoothly over $M$, hence so does its push-forward.
\qed

We next turn to the hermitian hessian itself, related to
the 2-form $\sigma_{\mlp}$ by the formula
$$g_{\mlp}(x,y)=\sigma_{\mlp}(x,Jy).$$  The following result is
part of a classical lemma usually attributed to Oka and Lelong.

\begin{prop}
\label{prop-oka}
If $\psi$ is any defining function for the strongly pseudoconvex
boundary $M$ of $X$, then $\mlp$ is strongly plurisubharmonic (i.e.
$g_{\mlp}$ is positive definite) on the
complement of $M$ in some neighborhood of $M$ in $X$.
\end{prop}

\pf
The quadratic
form $q(x)=g_{\mlp}(x,\overline{x})$ is expressed in terms of the
components of $x$ in our special basis by
$$q(x)=-(1/\psi)( a\gamma^0(x)\gambar^0(x) - b_k
\gamma^0(x)\gambar^k(x) - \overline{b}_k\gamma^k(x)\gambar^0(x) -
c_{jk}\gamma^j(x)\gambar^k(x))$$
$$ + (1/\psi^2)\gamma^0(x)\gambar^0(x).$$
The corresponding hermitian matrix is the positive function
$1/\psi^2$ times the $(1+n)\times(1+n)$ block
matrix with $1-\psi a$ in the upper left hand corner, $-\psi c_{jk}$ in
the lower right block, and $-\psi b_k$ and its adjoint
in the off diagonal row and column.  Since $c_{jk}$ is the matrix of
the Levi form, it is negative definite when $M$ is strongly
pseudoconvex.  By the Sylvester criterion (hermitian version),
the entire matrix will be positive definite if its determinant
is positive.  Expanding this determinant (if we ignore the overall factor
$1/\psi^2$) in minors of the top row, we obtain $\psi^n$ times the determinant of
$-c_{jk}$ plus terms divisible by $\psi^{n+1}$.  Sufficiently close to the
boundary, the sum must be positive.
\qed

\begin{prop}
\label{prop-quasi}
If (with notation as above) the function $\mlp$ is strictly
plurisubharmonic throughout the interior of $X$,
then the K\"ahler metric $g_\mlp$
on the interior of $X$ extends to a para-K\"ahler Lie algebroid
structure on  $E_{M,X}$.
\end{prop}
\pf
We have only to show that the summands in
the splitting of Example \ref{ex-splitting} are lagrangian with
respect to the symplectic form.  But this follows immediately by
continuity from the corresponding fact on the interior.
\qed
The nondegenerate pairing $\omega E^{1,0}_{M,X}\times
E^{0,1}_{M,X}\to\complex$
is given in our basis of sections by
$$\beta=(1-\psi a)\theta^0\otimes\theta'^0-b_k\theta^0\otimes \theta'^k
+ \overline{b}_k\psi\theta^k\otimes\theta'^0
-c_{jk}\theta^j\otimes\theta'^k.$$

\section{Quantization}
Deformation quantization on closed K\"ahler manifolds can be
accomplished as a by-product of Berezin-Toeplitz quantization.  See,
for instance, \cite{bo-me-sch:quantization} \cite{gu:star} and references
therein.\footnote{This work applies only to the case where the
  symplectic structure is integral, but Melrose \cite{me:star} has
  shown how to extend the method to the non-integral case.}
On the other hand, Karabegov \cite{ka:separation} studied special
formal deformation quantization adapted to the K\"ahler structure,
and he and Schlichenmaier \cite{ka-sch:identification} linked the
two approaches.

In this section, we will show how to extend the work cited above
to the case of K\"ahler manifolds with pseudoconvex boundary, using
para-K\"ahler Lie algebroids.

\subsection{Definitions}

We start with the ``complexification'' of a basic definition of Nest and Tsygan
\cite{ne-ts:quantization} (suggested already by Xu \cite{xu:quantum}).

\begin{dfn}
\label{dfn-estar}Let $(E, \omega)$ be a complex symplectic Lie
algebroid over $X$. An $E$-star product on $X$ is a formal series of
$E$-bidifferential operators
\[
B=1\otimes 1+\sum_{k\geq 1}(i\hbar)^k B_k,~B_k\in ^{E}\!\!
Op(X)\otimes ^{E}\!\!Op(X),
\]
which is associative in the sense that the Gerstenhaber bracket
$[B,B]$ is equal to zero, and for which the antisymmetrization of
$B_1$ is the $E$-bivector field $\pi$ inverse to $\omega$.
\end{dfn}

Pushing $B$ forward by the anchor of $E$, we obtain a formal
series of ($T_\complex X$)-bidifferential operators which gives a
star-product $\star$ for the Poisson structure which is the
push-forward of $\pi$.

Now we extend to the para-K\"ahler case the notion of quantization
with separation of variables.

\begin{dfn}
\label{dfn-bipolar}
If $E$ is a para-K\"ahler Lie algebroid, we call an $E$-star
product {\bf bipolarized} if the bidifferential operators $B_k$
all belong to ${^{E^{1,0}}\! Op(X)}\otimes {^{E^{0,1}}\!Op(X)}$.
\end{dfn}

A bipolarized $E$-star product has the property that $f\star g = fg$
whenever $f$ is an antiholomorphic function or $g$ is holomorphic.
When the anchor of $E$ is injective, even on a dense subset of $X$,
this property implies that the star product is bipolarized.  In the
K\"ahler case, this means that, after changing the sign of the complex
structure (or replacing the product by its opposite),
we are dealing with a star product with separation of variables in the
sense of \cite{ka:separation}.

\subsection{Bipolarized star products}
Nest and Tsygan \cite{ne-ts:quantization} showed that the quantization
method of Fedosov \cite{fe:simple} extends immediately to (real)
symplectic Lie algebroids to produce ``Weyl-type'' star products.
Their extension works for complex Lie
algebroids as well.  On the other hand, Neumaier
\cite{neu:quantization} showed that, when one starts the Fedosov
construction with a ``bipolarized'' symplectic connection (such as the
Levi-Civita connection for a pseudo-K\"ahler manifold) and a
``(anti)Wick-type'' bipolarized product on the tangent spaces, the resulting
star-product is bipolarized.  In this section, we combine the two
constructions above to obtain bipolarized $E$-star products on
para-K\"ahler Lie algebroids.

The idea of the construction of the $E$-star products can be
summarized as follows.

In this paragraph, we ``assume" that our para-K\"ahler Lie
algebroid (as a complex Lie algebroid) can be ``integrated" to a
$s$-connected ``groupoid" $G\rightrightarrows X$, which is a
``mystery" to the authors. The sections of our Lie algebroid can
be viewed as left $G$ invariant vector fields along the $s-$fibers
of $G$. The symplectic Lie algebroid structure defines a
$G-$invariant para-K\"ahler structure on each $s-$fiber.
Therefore, each $s$-fiber of $G$ becomes a para-K\"ahler manifold
and in particular a symplectic manifold. $G$ is canonically
foliated by the $s-$fibers with equal dimensions, and therefore
becomes a regular Poisson manifold, and $X$ is a complete
transversal to this foliation. A para-K\"ahler connection on $X$
can be lifted to a $G$-invariant symplectic connection on $G$.
Given a ``symplectic connection" on a regular Poisson manifold, we
can use Fedosov's construction (Neumaier's construction in
\cite{neu:quantization}) to obtain a star product on $G$. Since
the para-K\"ahler form and the para-K\"ahler connection are both
$G$-invariant, this product we obtain on $G$ is also
$G-$invariant, and therefore can be expressed by a $G$-invariant
bidifferential operator on $G$, which is a Lie algebroid
bidifferential operator and actually bipolarized if the
characteristic form is bipolarized.

The construction described in the above has one flaw that we have
assumed a complex Lie algebroid can be integrated. We can pass by
this issue by working with the $E-$jets introduced in Definition
\ref{dfn-ejets}, which can be viewed as the infinity jets of
smooth functions on $G$ restricted to $X$. The Grothendieck
connection defines a natural lift of the Lie algebroid action to
the infinity jets. Therefore, we can construct a bipolarized star
product by working with $E-$jets.

We begin our construction with the fibrewise anti-Wick product.
Let $E$ be a para-K\"ahler Lie algebroid over $X$.  Then each
fibre $E_x$ has a natural translation-invariant para-K\a"hler
structure given by the symplectic form and the lagrangian
subspaces $E^{1,0}_x$ and $E^{0,1}_x$.  Given a basis $\xi_0,
\ldots, \xi_{n}\in E^{1,0}, \xi'_0, \ldots, \xi'_n\in E^{0,1}$ and
the dual basis $\theta^0, \ldots, \theta^n \in {E^{0,1}}^*,
{\theta'}^0,\ldots, {\theta'}^n \in {E^{1,0}}^*$, $\omega$ is
expressed as $\omega_{ij}\theta^i \wedge {\theta'}^j$.

We define a bipolarized star product on the
algebra $$W_x\stackrel{def}{=}\complex[\eta^0,
\ldots, \eta^n, \eta'^0, \ldots, \eta'^n][[\hbar]]$$
of $\complex [[\hbar]]$-valued polynomial functions on $E_x$
by:
\begin{equation}
\label{eq-dfn-star} f\ast
g\stackrel{def}{=}\exp(-\frac{i\hbar}{2}\pi^{ij}\frac{\partial}{\partial
\eta'^i}\otimes \frac{\partial}{\partial \eta^j})(f\otimes g),
\end{equation}
where $(\pi^{ij})$ is the inverse matrix to $(\omega_{ij})$.
Taking a union of the algebras $W_x, x\in X$, we obtain a formal
anti-Wick algebra bundle $\calw$.

In the following, we adapt Fedosov's construction of star product
symplectic manifolds to our situation.

The para-K\"ahler connection introduced in Section
\ref{sec:para-conn} naturally lifts to a connection denoted  $\nabla^{lc}$
on the anti-Wick algebra bundle $\calw$.

\begin{dfn}
\label{dfn-flatconn}A Fedosov connection on $\calw$ is a flat connection
$D$ on $\calw$ of the form $D=\nabla^{lc} +A$, with $A\in \Omega^1(M,
End(\calw))$ and $D^2(a)=\frac{i}{\hbar}[\Omega, a]=0, \forall
a\in \Gamma(\calw)$.  $\Omega$ is a $\complex [[\hbar]]$-valued 2-form and usually
called the Weyl curvature of the connection $D$.
\end{dfn}

The following theorem is an extension of Fedosov's
Theorem on symplectic manifolds to para-K\"ahler manifolds..

\begin{thm}
\label{thm-flatconn} Let $\mu$ be an element of $
-\omega+\hbar\,{^{E^{1,0}}\Omega ^1(X, \complex})\wedge
{^{E^{0,1}}\Omega^1(X, \complex)}[[\hbar]]$ such that $d\mu=0$.
Then exists a $End(\calw)$-valued $E$-form $A_\mu$ on $X$ such
that $\nabla _{\mu}=\nabla^{lc}+A_{\mu}$ defines a Fedosov
connection on $\calw$ with
$\nabla_{\mu}A_{\mu}+\frac{1}{2}[A_{\mu}, A_{\mu}]=\mu$. The
complex $({^E \Omega (X, \calw)}, \nabla_{\mu})$ and
$({^E\Omega(X, {^EJets}\otimes \calw)}, \nabla_G+\nabla_{\mu})$
are acyclic in the positive dimension.
\end{thm}
$\pf$ The proof of this theorem is an application of Fedosov's
iteration method. The following construction is a generalization
of Theorem 3.1 in \cite{neu:quantization}. We outline the main
steps in the following, and the omit the detail check.

Notice that the anti-Wick algebra is naturally graded by the sum
of the degree of the polynomial and the power of $\hbar$. And
therefore, the algebra ${^E\Omega(X,\calw)}$ is graded by the
total degree, which is denoted by ``deg''. Following Fedosov
\cite{fe:simple}, we introduce operations $\delta$ and
$\delta^{-1}$ on ${^E\Omega(X,\calw)}$, as follows:
\[
\begin{split}
\delta(a)&=\sum_{i}(\theta^i\wedge\frac{\partial}{\partial
\eta^i}a+{\theta'}^i\wedge\frac{\partial}{\partial {\eta'}^i}a)\\
\delta^{-1}(a)&=\frac{1}{{\rm
deg}(a)}\sum_{i}(\eta^i\cdot\iota_{\theta^i}a+{\eta'}^i\cdot\iota_{{\theta'}^i}a),
\end{split}
\]
for homogeneous $a\in {^E\Omega(X, \calw)}$ with positive degree.
We look for $A_{\mu}$ of the form
$-\delta+\frac{i}{\hbar}[r_{\mu}, \cdot]$, where $r_{\mu}$ is in
${^E\Omega^1}(X, \calw)$ with total degree 5.

According to $\nabla_{\mu}A_{\mu}+\frac{1}{2}[A_{\mu},
A_{\mu}]=\mu$, $r_{\mu}$ is the unique solution of the following
equations
\begin{equation}
\label{eq:fe-conn}\delta
r_{\mu}=-\omega+\nabla^{lc}r_{\mu}+\frac{i}{\hbar}r_{\mu}\ast
r_{\mu}+R^{lc}-\mu,\  \delta^{-1}(r_{\mu})=0,
\end{equation}
where $R^{lc}$ is the curvature of $\nabla^{lc}$, which is in
${E^{1,0}}^*\wedge {E^{0,1}}^*$.

The solution of Equation (\ref{eq:fe-conn}) can be obtained by the
following iteration.
\[
r_{\mu}=\delta^{-1}(r_0)+\delta^{-1}(\nabla^{lc}r_{\mu}+\frac{i}{\hbar}r_{\mu}\ast
r_{\mu}),
\]
where $r_0=-\omega+R^{lc}-\mu$. \qed

By the flat connection $\nabla_{\mu}$ constructed in Theorem
\ref{thm-flatconn}, we have isomorphisms
\begin{equation}
\label{eq-qis-fun} \tau:C^{\infty}(X)[[\hbar]]\rightarrow
Ker(\nabla_{\mu}|_{^E\Omega^0(X, \calw)}),
\end{equation}
and
\begin{equation}
\label{eq-qis-jets} ({^E \calj(X)})\rightarrow
Ker((\nabla_G+\nabla_{\mu})|_{^E\Omega^0(X, ^EJets\otimes
^E\calw)}).
\end{equation}

By the isomorphism (\ref{eq-qis-jets}), the fiberwise product on
$^E\Omega^0(X, ^EJets\otimes ^EW)$ defines an associative product
$\ast$ on the space of $E$-jets. By the duality between the space
of $E$-jets and $^E Op(X)$, the product $\ast$ defines a left
$C^{\infty}(X)$-module map $\chi:{^EOp(X)}\to
{^EOp(X)}\otimes_{C^{\infty}(X)}{^EOp(X)}$, i.e. $l_1\ast
l_2(D)=l_1\otimes l_2(\chi(D))$. In particular, evaluating $\chi$
on the constant function $1$ in ${^EOp(X)}$, we obtain an
$^E-$bidifferential operator $\chi(1)$.

By arguments as in Theorem 4.3 of \cite{ne-ts:quantization}, we
have that $\chi(1)\in {^EOp(X)}\otimes{^EOp(X)}$ defines an
associative $E-$star product.

\begin{lemma}
\label{lem:polar}
\item For any $E-$jet $f$ constant along $E^{1,0}$ in an open set
$U$ of $X$, $\chi(1)(f, g)|_U=f\ast g|_U=fg|_U, \forall g\in
{^E\calj(X)}$.
\item For any $E-$jet $g$ constant along $E^{0,1}$ in an open set
$V$ of $X$, $\chi(1)(f, g)|_V=f\ast g|_V=fg|_V, \forall f\in
{^E\calj(X)}$.
\end{lemma}
$\pf$ When $E^{1,0}$ and $E^{0,1}$ are  the holomorphic and anti
holomorphic tangent bundle of a  K\"ahler manifold, Lemma
\ref{lem:polar} is Proposition 4.4 of \cite{neu:quantization}. The
present case is a small generalization thereof.  In the following,
we sketch the proof.

We study the isomorphism (\ref{eq-qis-fun}) in more detail. Let
$\hat{f}$ denote the lifting of an $E-$jet $f$ to a flat section
in $^E\Omega^0(X, ^EJets\otimes ^E\calw)$. $\hat{f}$ is the only
solution of the differential equation $\nabla_{\mu}\hat{f}=0$ with
$\hat{f}|_{\eta^i={\eta' }^i =0}=f$. We write
$\hat{f}=\sum_{k}\hat{f}^{(k)}$ into a direct sum of its
homogeneous $degree=k$ component. $\hat{f}^{(k)}$ can be
constructed through the following iteration,
\begin{equation}
\label{eq:lift-f}
\begin{split}
\hat{f}^0&=f,\\
\hat{f}^{(p)}=\delta^{-1}\big( \nabla^{lc}
\hat{f}^{(p)}&+\frac{i}{\hbar}\sum_{k=0}^{p-1}[r_{\mu}^{(k+2)},
\hat{f}^{(p-k)}]\big),
\end{split}
\end{equation}
where $r_{\mu}^{(k+2)}$ is the degree $k+2$ component of
$r_{\mu}$.

The product $\ast$ on ${^E\calj(X)}$ is then defined to be $f\ast
g\stackrel{def}{=}\sigma(\hat{f}\circ \hat{g})$, where
$\sigma:{^E\Omega}^0(X, ^EJets\otimes ^E\calw) \to {\Omega}
^{*}(X, {^E Jets})$ by setting all $\eta^i, {\eta'}^i$ equal 0 and
$\circ$ is the fiberwise multiplication in $^E\Omega^0(X,
^EJets\otimes ^E\calw)$.

To prove Lemma \ref{lem:polar}, we look at each fiber $W$ of the
bundle $\calw$. In Equation (\ref{eq-dfn-star}), we notice that
the product $\ast$ in the first component only involves $\eta'$,
and in the second component only involves $\eta$. Due to this
property of $\ast$, we introduce on $W$ a projection
$\tau^{1,0}:W\to W$ onto the component of $\complex({\eta'}^0,
\cdots, {\eta'}^n)[[\hbar]]$, the $(1,0)$ component, and a
projection $\tau^{0,1}:W\to W$ onto the component of
$\complex({\eta}^0, \cdots, \eta^n)[[\hbar]]$, the $(0,1)$
component. $\tau^{1,0}$, and $\tau^{0,1}$ lifts onto $\calw$
naturally, and denote the composition of $\tau^{1,0}$ and
$\hat{f}$ by $\hat{f}^{1,0}$ and similarly the composition of
$\tau^{0,1}$ and $\hat{g}$ by $\hat{g}^{0,1}$.

We observe that $f\ast g$ can be rewritten as
$\sigma(\tau^{1,0}(\hat{f})\circ \tau^{0,1}(\hat{g}))$, and
therefore it is enough to construct $\hat{f}^{1,0}$ and
$\hat{g}^{0,1}$ to compute $f\ast g$. Following the idea of
Proposition 4.2, \cite{neu:quantization}, we can restrict to
$E^{1,0}$ and $E^{0,1}$ to calculate $\hat{f}^{1,0}$ and
$\hat{g}^{0,1}$. They can be constructed by the following
iterations analogous to (\ref{eq:lift-f}).
\[
\begin{split}
\hat{f}^{1,0}&=f+\delta^{-1}_{1,0}(\nabla^{lc}_{1,0}\hat{f}^{1,0}+\frac{1}{\hbar}\tau^{1,0}(\hat{f}^{1,0}\circ
r_{\mu}^{1,0})),\\
\hat{g}^{0,1}&=g+\delta^{-1}_{0,1}(\nabla^{lc}_{0,1}\hat{g}^{0,1}+\frac{1}{\hbar}\tau^{0,1}(\hat{g}^{0,1}\circ
r_{\mu}^{0,1})),
\end{split}
\]
where $\nabla^{lc}_{1,0}, \nabla^{lc}_{0,1}$, and $r_{\mu}^{1,0},
r_{\mu}^{0,1}$ are the restrictions of $\nabla$ and $r_{\mu}$ to
their components in $E^{1,0}$ and $E^{0,1}$.

With the formulas above, it is straightforward to check that if
$f$ is constant along $E^{1,0}$, the $\hat{f}^{1,0}=f$ and $f\star
g=fg$, and the same holds for the $E^{0,1}-$component. $\Box$

\begin{prop}
\label{prop-quantization} The isomorphisms $(\ref{eq-qis-fun})$
and $(\ref{eq-qis-jets})$ define a bipolarized $E$-defor\-mation of
a para-K\"ahler Lie algebroid.
\end{prop}
$\pf$ We need to show that $\chi(1)\in {^{E^{1,0}}Op(X)}\otimes
{^{E^{0,1}}Op(X)}$, but this is implied by Lemma \ref{lem:polar}.
\qed

\subsection{Strictly pseudoconvex boundary of a complex manifold}
Corollary \ref{cor-poisson} shows that the para-K\"ahler Lie
algebroid associated to the boundary of a complex domain has a
polarized $E$-star product. On the other hand, in
\cite{en:berezin}, Engli\v s constructed a differential star
product on a bounded pseudoconvex domain by Berezin-Toeplitz
quantization. In this subsection, we identify Engli\v s' star
product with what we constructed in Corollary \ref{cor-poisson}.

Engli\v s' Berezin-Toeplitz quantization is an example of a
deformation quantization with separation of variables. In Theorem
5.9 \cite{ka-sch:identification}, Karabegov and Schlichenmaier
identified the opposite of the Berezin-Toeplitz star product with
a star product with separation of variables whose Karabegov form
is equal to $-\frac{1}{i\hbar}\omega+\omega_{can}$ and
characteristic form is
$\frac{1}{i\hbar}\omega-\frac{\omega_{can}}{2i}$, where
$\omega_{can}$ is the curvature form of the canonical line bundle.

\begin{prop}
\label{prop-canonical} Near the boundary $M$ of a pseudoconvex
domain $X$, the ca\-nonical form  $\omega_{can}$ pulls back to a smooth section of
${E^{1,0}}^*_{M, X}\otimes {E^{0,1}}^*_{M, X}$.
\end{prop}
\pf
 According to the definition given in
$\cite{ka-sch:identification}$, $\omega_{can}$ is equal to $-i
\partial \overline{\partial} \nu$, where $i^{n+1}e^{\nu}dzd\bar{z}$
 is
the symplectic volume.  (We write $dz$ for $dz^1\cdots dz^{n+1}$.)
To calculate $\omega_{can}$, we follow the
calculation in Proposition \ref{prop-nondegen}, where the volume
form is shown to be
\[
\left(-\frac{\del\delbar\psi}{\psi} +
\frac{\del\psi\wedge\delbar\psi}{\psi^2}\right)^{n+1}
=\left(-\frac{\del\delbar\psi}{\psi}\right)^{n+1}+
\left(-\frac{\del\delbar\psi}{\psi}\right)^{n}\wedge
\frac{\del\psi\wedge\delbar\psi}{\psi^2}.
\]
This form can be written as
$\psi^{-n-2}\frac{det(\psi\partial_i\delbar_j
\psi-\del_i\psi\delbar_j \psi)}{\psi ^{n}}dzd\bar{z}$. It has been
shown in Theorem \ref{thm-algebroid} that
$\nu_0=\frac{det(\psi\partial_i\delbar_j \psi-\del_i\psi\delbar_j
\psi)}{\psi ^{n}}$ is nonzero and smooth up to the boundary.
Therefore, $\nu$ is equal to
$\log(\psi^{-n-2}\nu_0)=-(n+2)\log(\psi)+\log(\nu_0)$, and
$\del\delbar\nu$ is
\[
-(n+2)\del\delbar\log (\psi)+\del\delbar(\log (\nu_0)).
\]

To show that $\del\delbar \nu$ pulls back to a smooth section of
${E^{1,0}}^*_{M, X}\wedge{E^{0,1}}^*_{M, X}$, we may check
the two terms separately.
\begin{enumerate}
\item $-(n+2)\del\delbar\log (\psi)$. By Theorem \ref{thm-algebroid}, $\del\delbar \log(\psi)$ pulls
back to be a smooth section of ${E^{1,0}}^*_{M,
X}\wedge{E^{0,1}}^*_{M, X}$.
\item $\del\delbar(\log (\nu_0))$. Since $\nu_0$ is nonzero near the boundary, $\del\delbar
\log(\nu_0)$ is smooth up to the boundary, and therefore expressed
in terms of $\gamma^i$ and $\gambar^i$, which are smooth sections $\wedge^2
E^*_{M, X}$.  Finally, the bipolarization of $\del\delbar(\log(\nu_0))$ is
obvious. \qed
\end{enumerate}

Since $\omega_{can}\in {^{E^{1,0}}\Omega ^1(X, \complex})\wedge
{^{E^{0,1}}\Omega^1(X, \complex)}$, we obtain a bipolarized
$E$-star product $\star$ on $X$ constructed from Proposition
\ref{prop-quantization} with the Weyl curvature $\mu$ equal to
$\mu=-\omega+i\hbar\omega_{can}$.
\begin{prop}
\label{prop:karabegov-form} When restricted to the interior of
$X$, the Karabegov form of $\star$ is equal to
$-\frac{1}{i\hbar}\omega+\omega_{can}$.
\end{prop}
$\pf$ When restricted to the interior of $X$, $E_{M,X}$ coincides
with the tangent bundle of the K\"ahler manifold $X$. By the
locality of our construction of star products, we have that in the
interior of $X$, our quantization of a para-K\"ahler Lie algebroid
$E_{M,X}$ coincides with the quantization of the K\"ahler
manifold. This allows us to use Theorem 6.7 and Deduction 6.9 in
\cite{neu:quantization} to obtain the characteristic form. $\Box$

Now we are ready to state the following theorem.
\begin{thm}
\label{thm:smooth-boundary} The opposite Berezin-Toeplitz star
product $\star'_{BT}$ near a strictly pseudoconvex boundary is
equal to the $E$-star product $\star_E$ constructed in Proposition
\ref{prop-quantization}, and therefore is smooth up to the
boundary of $X$.
\end{thm}
\begin{rmk}
{\em
Smoothness up to the boundary for Engli\v s' Berezin-Toeplitz star product 
can be directly derived from Corollary
\ref{cor-para}.  (See Remark \ref{rmk-smoothness}.)
But Theorem \ref{thm:smooth-boundary} sets up a nice connection between the
Berezin-Toeplitz quantization and deformation quantization.
}
\end{rmk}
$\pf$ We have known that both the $E-$star product $\star_E$ and
the opposite Berezin-Toeplitz star product $\star'_{BT}$ in the
interior of $X$ are star products with separation of variables
defined by Karabegov. Karabegov in Theorem 2, \cite{ka:separation}
showed that star products with separation of variables on a
K\"ahler manifold $X$ are one to one correspondent to the set of
power series in $\hbar$ of closed $(1,1)$ forms on $X$,
i.e.$Z_{dR}^2(X, \complex)^{(1,1)}[[\hbar]]$, which is called
Karabegov form.

We have shown in Proposition \ref{prop:karabegov-form}, the
Karabegov form of the $E-$star product $\star_E$ is equal to
$-\frac{1}{i\hbar}\omega+\omega_{can}$. A similar calculation as
Karabegov and Schlichenmaier did on compact K\"ahler manifolds
computes that the Karabegov form of Engli\v s' opposite
Berezin-Toeplitz star product $\star'_{BT}$ is also equal to
$-\frac{1}{i\hbar}\omega+\omega_{can}$. Therefore, we have
\[
\star'_{BT}=\star_E.
\]
And since $\star_E$ is smooth up to the boundary of $X$,
$\star'_{BT}$ is also smooth up to the boundary, which was already
proved in Corollary \ref{cor-para}.

We end the proof by explaining the computation of the Karabegov
form of Engli\v s' opposite Berezin-Toeplitz star product.

Karabegov and Schlichenmaier's calculation in
\cite{ka-sch:identification} requires the K\"ahler manifold under
consideration to be compact.  In \cite{ka-sch:identification}, 
they use the compactness assumption 
to prove that the Berezin-Toeplitz quantization defines
a local differential product.  Their Propositions 5.1
and 5.2 give an asymptotic expansion of the Berezin
transform and the twisted operator product, while Theorem 5.6
gives an asymptotic expansion
of the Bergman kernel.  All of the remaining calculations are local.
For the case of the boundary of a strictly pseudoconvex domain, the
analogs of the above asymptotic expansions have been established by
Engli\v s  in  Theorem 1-3 of  \cite{en:berezin}. Therefore, we can
still use the method of formal integrals as in
\cite{ka-sch:identification} to compute the Karabegov form of the
star product $\star'_{BT}$. It turns out be again equal to
$-\frac{1}{i\hbar}\omega+\omega_{can}$. $\Box$

\end{document}